\documentclass[a4paper,11pt]{article}
\usepackage{amssymb,rotating,graphics,epsfig,float,graphics}
\usepackage{amsmath,verbatim,a4wide}
\usepackage{theorem}
\usepackage{color}

\definecolor{blue}{rgb}{0,0,1}
\definecolor{red}{rgb}{1,0,0}
\definecolor{purple}{rgb}{1,0,1}

\long\def\red#1\endred{\textcolor{red}{#1}}
\long\def\blue#1\endblue{\textcolor{blue}{#1}}
\long\def\purple#1\endpurple{\textcolor{purple}{#1}}

\theoremstyle{break} 
{\theorembodyfont{\rmfamily} }
{\theorembodyfont{\rmfamily} }
{\theorembodyfont{\rmfamily} }
{\theorembodyfont{\rmfamily} }

\begin{document}

\author{Ferdinand Verhulst \\
Mathematisch Instituut,  Utrecht University \\
PO Box 80.010, 3508TA Utrecht, The Netherlands} 

\title{A conceptual model for\\ growth by Capital-Education investments} 

%\date{  }

\maketitle 

\begin{abstract} 

Economic growth depends on capital investments and on investments 
in education and innovation. The model introduced here will specifiy aggregate output as determined by 
aggregate supply of capital and education investment. 
After formulating and analysing such a model in section 2 we will  consider the 
effectiveness of 
education for the growth of the National Product. It turns out that small changes 
of the quality of education has a  considerable impact on economic growth. Secondly we consider 
the influence of chaotic fluctuations of capital investments caused by hype-cycles or erratic 
policies. In section 3 we introduce a continuous control 
on education investments depending on 
consumption. In this 3-dimensional macro-economic model it turns out that a tipping point exists 
where increase of consumption affecting the amount of 
education and innovation leads to decline of economic growth.
\end{abstract} 


JEL Class. H52, I25, O38; \\
MSC Class. 91B62
\section{Introduction}  
The development of
macro-economic models in the middle of the 20th century produced a great many interesting 
discussions. A seminal paper by Solow \cite{S56} in 1956 gave a critical assessment of the classical 
Harrod-Domar model while describing useful modifications and extensions. 
Phelps discussed in 1961 models for maximizing consumption in a growing economy, 
see \cite{P61} and section \ref{sec3}. 
The models are derived making assumptions on quantities that change in time like capital investments 
and quantities and relations that are independent of time or at least quasi-stationary like 
the effectivity of investments expressed by elasticity coefficients.   
In \cite{H76} Hicks distinguished in 1976 between economic theory that involves time and so-called `economic 
theory out of time'. This  distinction will play a part in our section \ref{sec2} on  a basic model 
for the time-evolution of the National Product. 
Somewhat later, in 1992, an interest arousing paper \cite{MRW} appeared with arguments that  
schooling, say education, is  closely tied to economic growth. This was discussed  again in  
\cite{BG01}.  
The report \cite{B11}  discusses in detail the various social aspects and interactions of workers 
with schooling and its impact on te National Product. It introduces an appropriate Cobb-Douglas production function with estimates of the 
essential parameters, it reviews cross-country studies of the elasticity coefficients of capital $K$ 
and human capital in the form of education and schooling of workers.  There is relative agreement 
about the elasticity coefficient of capital investments but more variation on the effectiveness of 
eeducation expressed by\ its eleaticity coefficient. 
In later publications, education as an enabler of economic growth has been discussed
 from many different points of view 
for instance regarding the quality of education, differences between developed and less 
developed countries, also as one of the causes of 
the fast growth of the economy of China. Recent and upcoming aspects are regularly updated in the 
OECD reports \cite{OECD} with headings Analyse by country, Explore data, Review education 
policies.\\
There are many political discussions about funding education, questioning whether increasing 
funds for schools 
improves education  or even to consider it as a waste of money. In \cite{K18} the arguments against education are formulated. This sounds extreme but there is probably a difference between 
a society where participation in education is completely state-funded without restrictions on the choices
 of schooling and a society where participation in education asks for significant financial 
 contributions from students and 
 where  pressure may exist on the choices of schooling.  These aspects (`learning something that 
 matters') are considered again in 
discussions on the `knowledge capital' of society and the quality of schooling in relation to 
economic growth, see \cite{HW15} and \cite{H08}. \\ 
The connection between skills and economic growth is one aspect, it is often classified under `Human 
Resources' but we will include innovation as a direct extension of education at universities and in 
industry. This will also take into account laboratories and material instruments as ingredient. 
The consequences of changes of government policies for the quality of education and for economic 
growth will be one of the aspects to discuss in the sequel.\\
Another aspect is the part played by timescales. Actions by Central Banks or government decisions 
to change VAT percentages have an immediate effect. Investments in education and innovation 
can take years to become effective. When we focus on government policies to finance education 
we will take care of the timescale issue by introducing continuous control as a new instrument.

\subsection{Set-up of the paper}
We will present a conceptual macro-economic model with basic variables capital $K$, available 
education and research activities $E$ with corresponding investment parameters $s_k, s_r$.  
Economic reality has a high grade of complexity and it is not easy to find out from large-scale 
models involving hundreds of variables and parameters what the result is from certain technical 
changes 
or the effect of changes by political decisions. \\
A conceptual model, see \cite{B11} or the classic text \cite{S48}  can, as in physics or 
engineering, elucidate the part played by 
certain key variables and parameters. A conceptual model also helps to study 
the impact of political choices  on the National Product by increaing or decreaseing the 
investments in education and research. 
Such a conceptual model neglects micro-economic management and many other aspects but 
it helps to clarify the overall development of the economy and the consequences of certain political choices. \\

We will use the following macro-economic concepts that are a function of time :
\begin{itemize}
\item National Income or National Product $Y$. This is the total amount of money spent per year 
by government, business and private persons on capital goods, education, other production 
means, goods  and services. 
\item Consumption $C$: the total amount of money spent per year by the population on goods 
and services.
\item Capital $K$: the total of physical production means valued in money.
\item Investment  $I_k$: the total amount of money spent per year on updating,  repairing and 
replacing capital goods. 
\item The total amount $E$ of education, expertise and research present in the population, 
including physical goods, valued in money (in a number of publications called 
`human capital'). 
\item Investment  $I_r$: the total amount of money spent per year on education and research. 
\end{itemize} 
The macro-economic quantities $Y, C, K, E$ are functions of time with a certain time unit; here the 
unit is a year but this is arbitrary. 
The control factor $p$  can be used by the government to determine the fraction of the National 
Product that can be spent on consumption, so $C=pY$.\\
The modeling will 
take place by considering conservation laws and empirical rules. \\ 
In section \ref{sec2} we will formulate a Cobb-Douglas production function $Y$ involving the quantities 
$ K$ and $E$. The elasticity parameters $\alpha, \beta$ and decay parameters $\delta_k, \delta_r$ 
will represent the effectiveness of education together with 
research and capital investments coefficients $s_r, s_k$. The phase-plane depicted 
in fig.~\ref{fig1} shows the dynamics of the basic model 
with  one stable equilibrium if $0 < \alpha + \beta <1$ . 
We expect that small social and political changes will influence capital investments. This is 
modeled in section \ref{sec2} by a chaotic timeseries. \\
In section \ref{sec3} we assume that the government wants to control the investments in education 
and research to guarantee the consumption $C$ as a fraction $p$ of the National Product $Y$. 
Because of the long timescale of these investments we choose a continuous control of the 
investment  parameter $s_r$.  The consequences of the choice of the fraction $p$ are clearly 
illustrated in fig.~\ref{fig4a}. \\
As discussed in \cite{B11}, \cite{HW15} and \cite{H08} the effectiveness of education may vary 
according to 
policy decisions. We will discuss the surprising consequences of a small decrease of effectiveness 
of education by varying the control coefficient $p$ in section \ref{sec3}. Economic growth turns out 
to be very sensitive to this control.

\section{The basic Capital-Education model} \label{sec2}
We will have the conservation law that National Income $Y$ equals the sum of Consumption $C$ and 
total investment, in our set-up $(I_k+I_r)$: 
\begin{equation}
 Y = C + I_k + I_r. 
 \end{equation} 
So both investments are fractions of National Income, we can express this as:
\begin{equation} \label{invest}
I_k = s_k Y, \, I_r = s_r Y;  0 \leq s_k \leq 1, 0< \delta \leq s_r \leq 1. 
\end{equation} 
We have $s_r \geq \delta$ with $\delta$ a small positive parameter as all mammals teach their 
young at least to forage. 
An empirical law is the {\em law of diminishing returns}, stating the observation that in 
general extra input of capital, education and innovation produces more growth in a sublinear way. 
A generalized Cobb-Douglas production 
function expresses this by assuming that the National Income  $Y$ is proportional to the schooling 
and expertise $E$ of the workforce and physical capital $K$ as:
\begin{equation} \label{CB}
Y =E^{\alpha}K^{\beta},\, 0 < \alpha, \beta < 1,
\end{equation} 
with $\alpha, \beta$ the elasticity coefficients. $E^{\alpha}$ is a productivity factor dependent on 
 the level of education and research  of the population and in its turn it is 
also dependent on the investment factor $s_r$. We will use in examples typical values of the 
elasticity coefficients following the cross-country survey in \cite{B11}; typical would be values 
near $\alpha =0.20,  \beta = 0.35$.\\

\begin{figure}[ht]  
\begin{center}
\resizebox{!}{8cm}{
\includegraphics{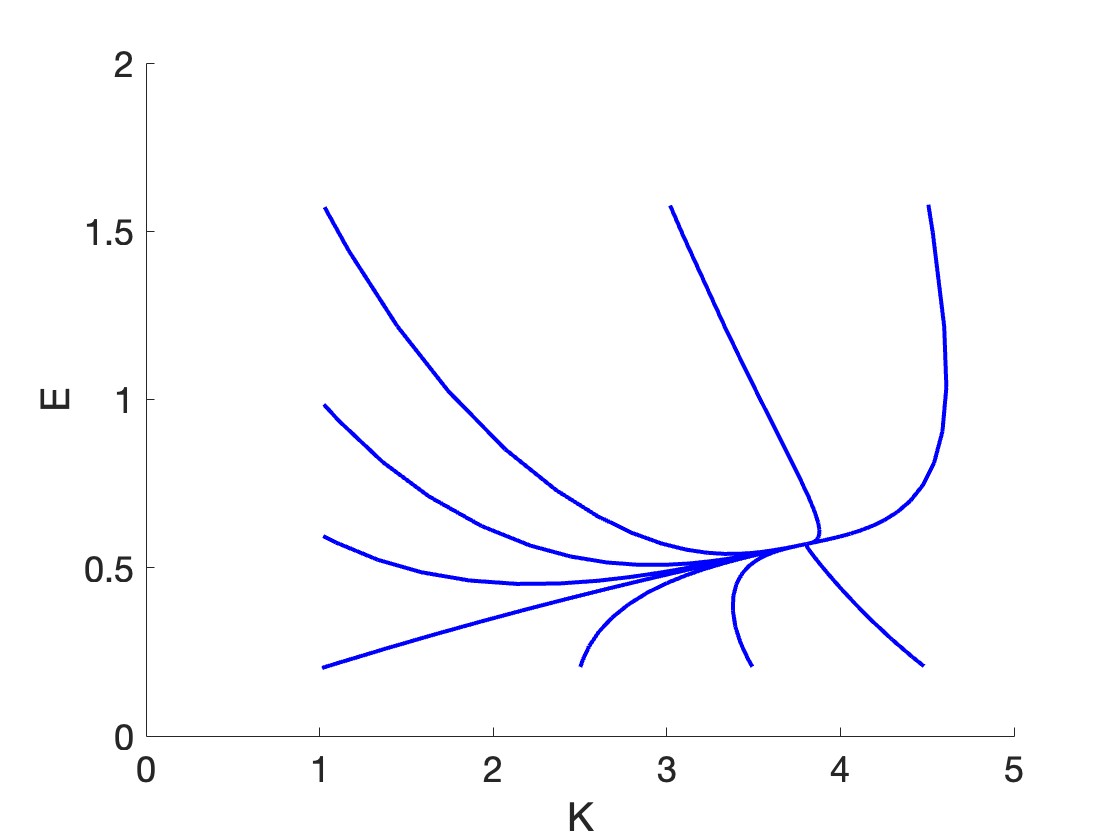}} 
\end{center}
\caption{Phase-plane dynamics of system \eqref{sys1}. Parameter values: 
$ s_r=0.1, \delta_r= 0.25, s_k=0.4, \delta_k=0.15,   \alpha = 0.2,  \beta=0.35$. 
Stable equilibrium at $(K, E)= (0.38, 0.57)$.
 \label{fig1}}
\end{figure} 
Of course there are many more aspects regarding the development of National Income but we intend 
in this note to consider the consequences of changes of investment in education, research and 
innovation. If 
necessary we can add small chaotic modulations of capital investments, this produces small 
modulations of the time series $Y(t)$.\\
The amount of capital $K$ increases by investment and decreases through wear and tear 
($\delta_k$-proportional to $K$). This leads to the equation
\begin{equation} \label{dKdt}
\frac{dK}{dt} = I_k - \delta_k K, \, \delta_k >0.
\end{equation}
Using eqs. \eqref{invest},   \eqref{CB} we find: 
\begin{equation} \label{dkdt}
\frac{dK}{dt} = s_k E^{\alpha}K^{\beta} - \delta_k  K.
\end{equation}
The parameters $s_k, \delta_k$ are semi-definite positive constants; it is natural that the parameters
 will change with time, but at least for some time we will 
assume that they are constant. For the productivity factor $E^{\alpha}$ we have in a similar way 
the equation:
\begin{equation} \label{prod}
\frac{dE}{dt}= s_rY - \delta_r E,  \delta_r>0,
\end{equation} 
with $\delta_r$ a positive parameter representing loss by obsolete and forgotten expertise. 
Using eqs. \eqref{invest}, \eqref{prod}, \eqref{dKdt} we find the system:
\begin{eqnarray} \begin{cases} \label{sys1} 
\frac{dK}{dt}&  = s_k E^{\alpha}K^{\beta} - \delta_k  K,\\
\frac{dE}{dt} & = s_r E^{\alpha}K^{\beta} - \delta_r E.
\end{cases} \end{eqnarray} 

If accidentally $\delta_k=\delta_r$, system \eqref{sys1} has an invariant manifold simply described by: 
\begin{equation} \label{figE1}
K= \frac{s_k}{s_r} E.
\end{equation} 
However, $\delta_k$ and $\delta_r$ are not related so we will pay no attention to this case. 
 We will consider the $(E, K)$ phaseplane. 
Apart from the trivial solution $(E, K)= (0, 0)$ we have one critical point $(E_0, K_0)$ given by:  
\begin{equation} \label{critpt}
E_0^{\alpha} = \frac{\delta_k}{s_k} K_0^{1- \beta},\,
 K^{\beta}_0=  \frac{\delta_r}{s_r} E_0^{1 - \alpha} . 
\end{equation} 

\begin{figure}[ht]  
\begin{center}
\resizebox{!}{5cm}{
\includegraphics{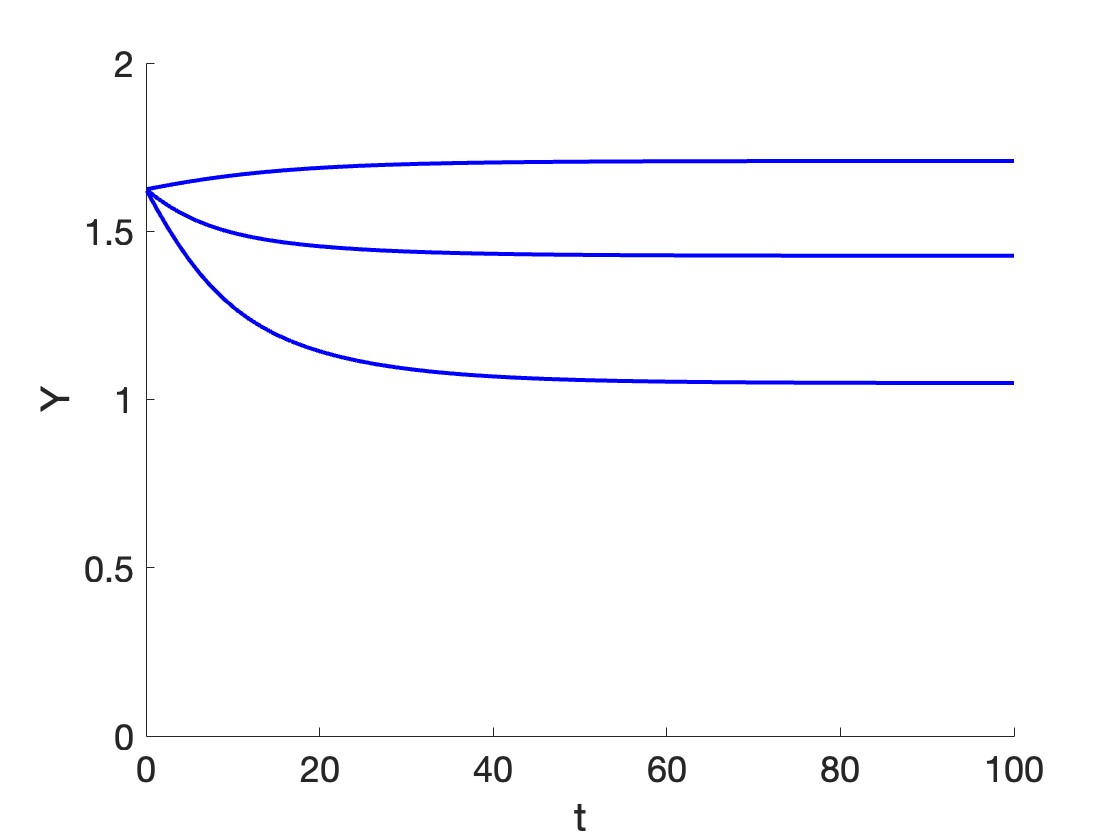}} 
\resizebox{!}{5cm}{
\includegraphics{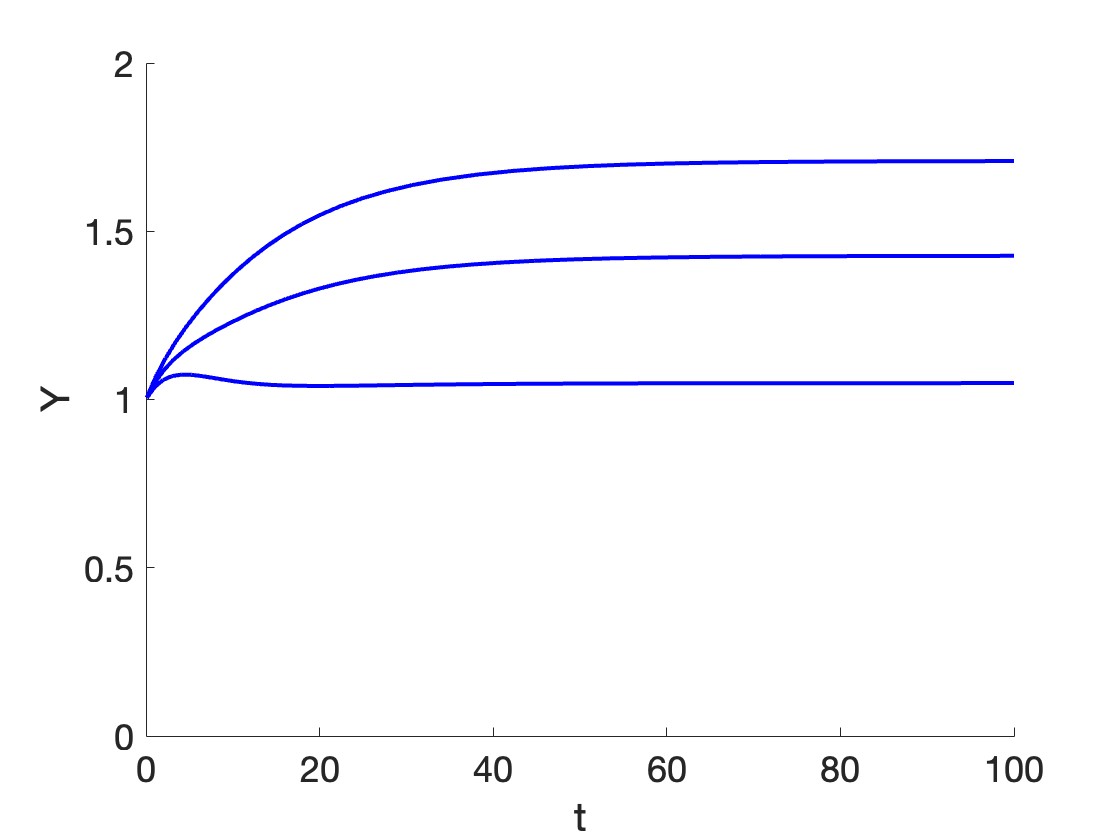}} 
\end{center}
\caption{ Time series of the National 
Product $Y(t)$  based on system \eqref{sys1} for initially $K(0)/ E(0)=4:1$ (left) and various  values of $s_r$; the parameter values are as for fig.~\ref{fig1}
and successively $s_r=0.05, 0.1, 0.15$. Decreasing growth is found at $s_r= 0.05, 0.1$, 
highest growth at $s_r= 0.15$,
Right the $Y(t)$ time-series with $K(0)=E(0)=1$ showing related but 
quantitatively different patterns. 
\label{fig2}}
\end{figure} 

This equilibrium (critical point) does not exist for arbitrary parameter values if $\alpha + \beta =1$; 
it would not be practical to impose this relation as the 2 elasticity coefficients are qualitatively 
very different.  Such a special 
choice of parameters produces in a 
simple way the AK model that allows for perpetual, exponential growth, assuming 
relatively large positive
savings and investment rates and $\delta_k, \delta_r$ sufficiently small; see also \cite{Ro86} and \cite{R91}. The case is called 
``structurally unstable'' in the mathematical theory of dynamical systems. This means that the 
case is exceptional in the sense that small perturbations will produce large quantitative en 
qualitative changes (see for the terminology \cite{FV25}). \\
Consider the stability of the economic equilibrium $(E_0, K_0)$ .
Linearising system \eqref{sys1} at the critical point $(E_0, K_0)$ gives a 2-dimensional matrix 
 (excluding the case  $\alpha + \beta =1$):
 \begin{equation} \label{lincrit}  \left( \begin{array}{cc}
(\beta  - 1) \delta_k & \alpha \frac{s_k}{s_r} \delta_r \\
\beta \frac{s_r}{s_k} \delta_k & (\alpha -1) \delta_r
\end{array}  \right) \end{equation}
with  trace representing the divergence of the flow near the critical point:
\begin{equation} \label{trace} 
(\alpha -1) \delta_r + (\beta -1) \delta_k <0.
\end{equation} 
 As the divergence is in the linearised system 
negative, the flow is locally contracting. We shall show that the corresponding economic 
equilibrium  is a stable node if $\alpha + \beta < 1$, see fig.~\ref{fig1}. 
The parameter values are suggested by the 
cross-country survey of \cite{B11} with typical investments ratio $s_k/s_r=4:1$.  \\

{\bf Stability analysis of the equilibrium $E_0, K_0$}.\\ 
As the divergence given by eq. \eqref{trace} is negative we have for the eigenvalues 
$\lambda_1, \lambda_2$ of matrix \eqref{lincrit} that $\lambda_1, + \lambda_2 <0$. 
From the characteristic equation we have  $\lambda_1 \lambda_2 = (1 - \alpha - \beta) 
\delta_r \delta_k$.
So one of the eigenvalues is positive if  
\begin{equation} 
\alpha + \beta >1.
\end{equation} 
In this case the equilibrium $(E_0, K_0)$ is unstable. Regarding he cross-country survey in \cite{B11} 
we consider these $\alpha, \beta$ values as less realistic. If we would have $\alpha + \beta > 1$ 
we would have the possibility of tipping points and  permanent growth of the National Product.\\

We assume from now on $ 0 <\alpha + \beta <1$.\\ 
In this case both eigenvalues $\lambda_1, \lambda_2$  are real and negative. 
Three typical orbits for the evolution of the National Product $Y(t)$ are shown in fig.~\ref{fig2}. 
In fig.~\ref{fig1} it is shown that the $E, K$ phase-orbits tend to a stable equilibrium, so as 
expected the National Product $Y(t)$ stabilises at a definite value. 
We use the parameters of fig.~\ref{fig1} in fig.~\ref{fig2} except for different values of investment 
coefficient $s_r$. Lowest growth of the National Product $Y(t)$ takes place at $s_r= 0.05$, slightly 
better for $s_r= 0.1$ and again increase if $s_r= 0.15$.  Starting with an equal initial ratio 
$K/E$ produces a remarkable growth of $Y(t)$ in fig.~\ref{fig2} (right). 

\subsection{The effectiveness of education}  \label{sec4b} 

\begin{figure}[ht]  
\begin{center}
\resizebox{!}{8cm}{
\includegraphics{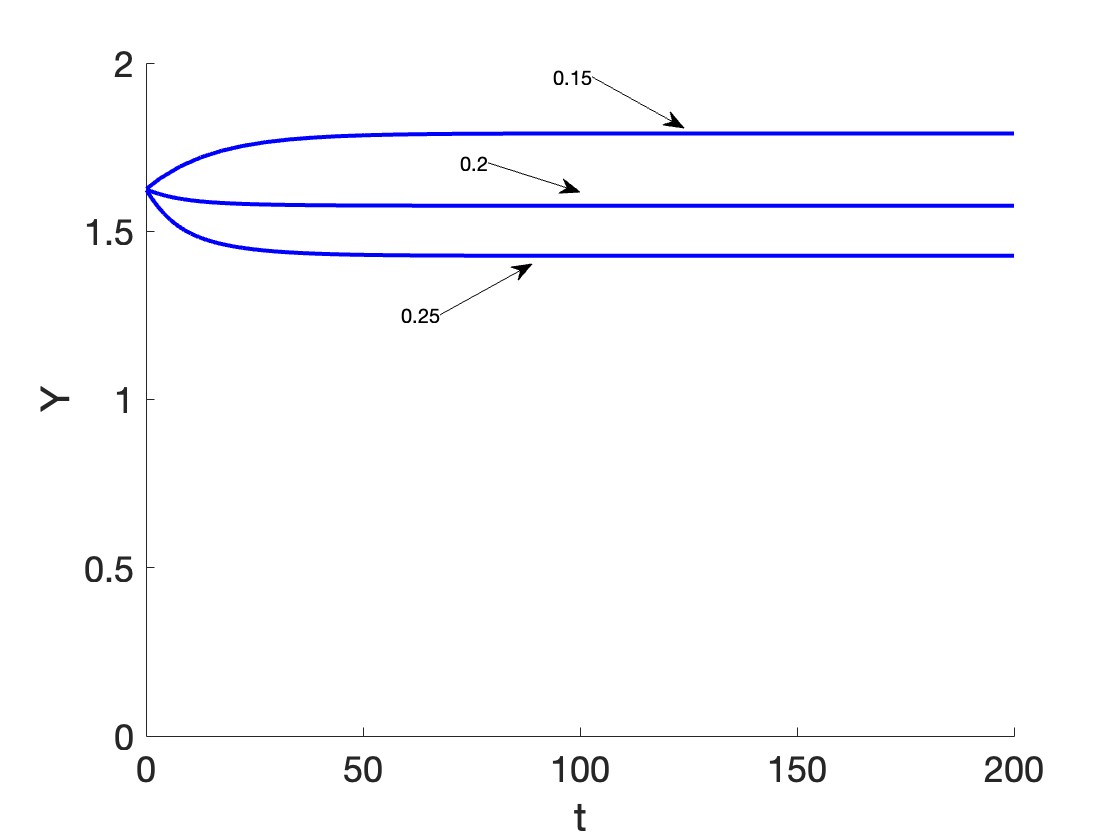}} 
\end{center}
\caption{National Product $Y(t)$ with parameter values 
$ s_r=0.1,  s_k=0.4, \delta_k=0.15,   \alpha = 0.2,  \beta=0.35$. 
The parameter $\delta_r$ will be varied to show the result of more effective education. 
The value $\delta_r =0.25$ corresponds with the time-series in the middle of fig.~\ref{fig2} (left). 
 Putting $\delta_r =0.2$ gives better growth results, putting $\delta_r =0.15$ even more.
 \label{fig2c}}
\end{figure} 
As discussed in \cite{H08}, \cite{HW15}, \cite{SR} and in numerous newspaper articles the 
investments in and the quality of 
education and research 
play an important part in economic growth. Regarding `quality' we include communication tools 
(languages), mathematics,  up-to-date tools from information theory and other modern topics, 
also specialized professional training by academic institutes and polytechnics. The ability to apply 
these tools in real-life situations is also essential. \\ 
An aspect of the discussion in more developed countries is also part-time work of highly qualified 
workers. Especially in teaching and medical 
professions the number of parttimers in developed countries has increased considerably. This leads 
clearly to less effectiveness as the cost to educate and  professionalise a person 
remains the same regardless of the number of hours worked when professionally active. \\ 

In fig.~\ref{fig2c} we show the results of more effective education by varying $\delta_r$. 
One could also 
use the elasticity coefficient $\alpha$ but then, because of the fractal exponent near zero, one has 
to choose the other parameters such that 
initially $dY/d \alpha >0$; with such a choice growth of $Y(t)$ with $\alpha$ is guaranteed. 
The lowest growth curve in fig.~\ref{fig2c} is taken from fig.~\ref{fig2}, decreasing $\delta_r$ 
produces more growth $Y(t)$. 
Methods to improve the effectivity of education are a topic of permanent discussion among 
teachers and educationalists.  The following National Product values were obtained after 200 
timesteps when varying $d_r$ and keeping the other parameters as in fig.~\ref{fig2c}: \\
\begin{equation} \begin{array}{ccccccc} \label{effectdr} 
d_r &  0.25 & 0.23 & 0.21 & 0.19 & 0.17 & 0.15 \\
Y & 1.43 & 1.48 & 1.54 & 1.61 & 1.68 & 1.79
\end{array} \end{equation}  
Increasing the effectivity of education has serious consequences for economic growth. 

\subsection{Chaotic fluctuations of investments} 
\begin{figure}[ht]  
\begin{center}
\resizebox{!}{5cm}{
\includegraphics{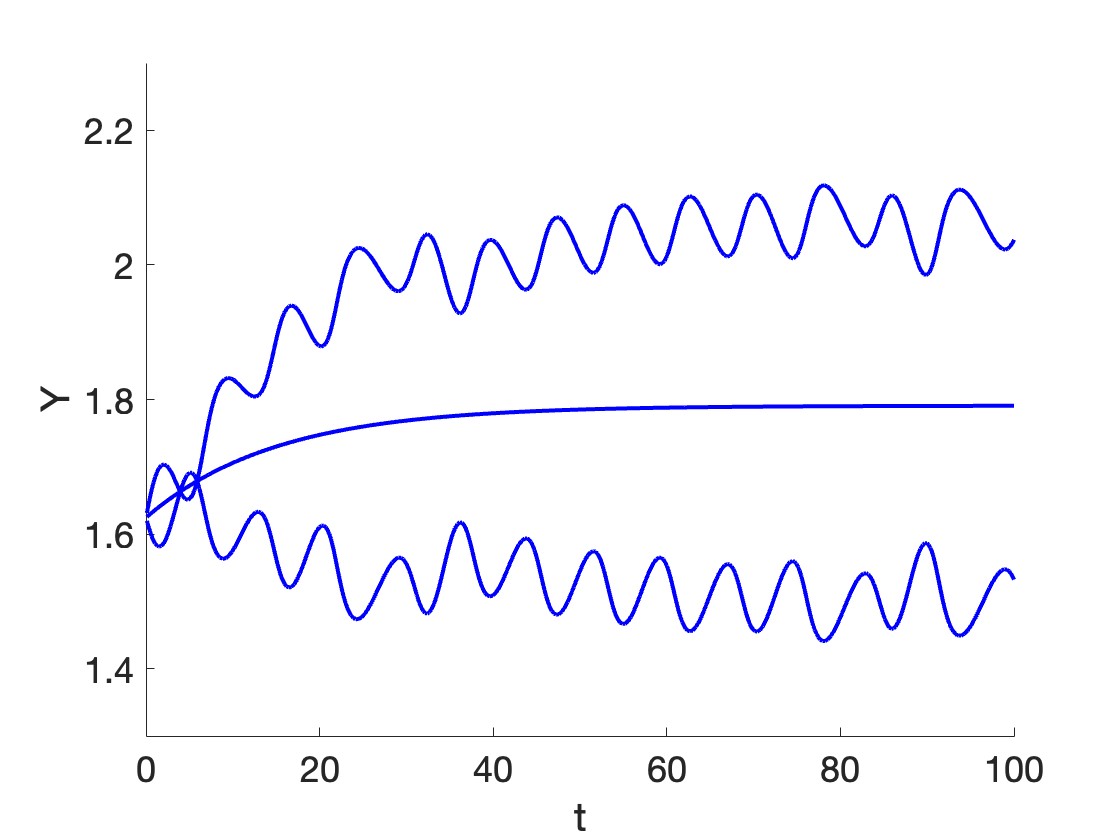}} 
\resizebox{!}{5cm}{
\includegraphics{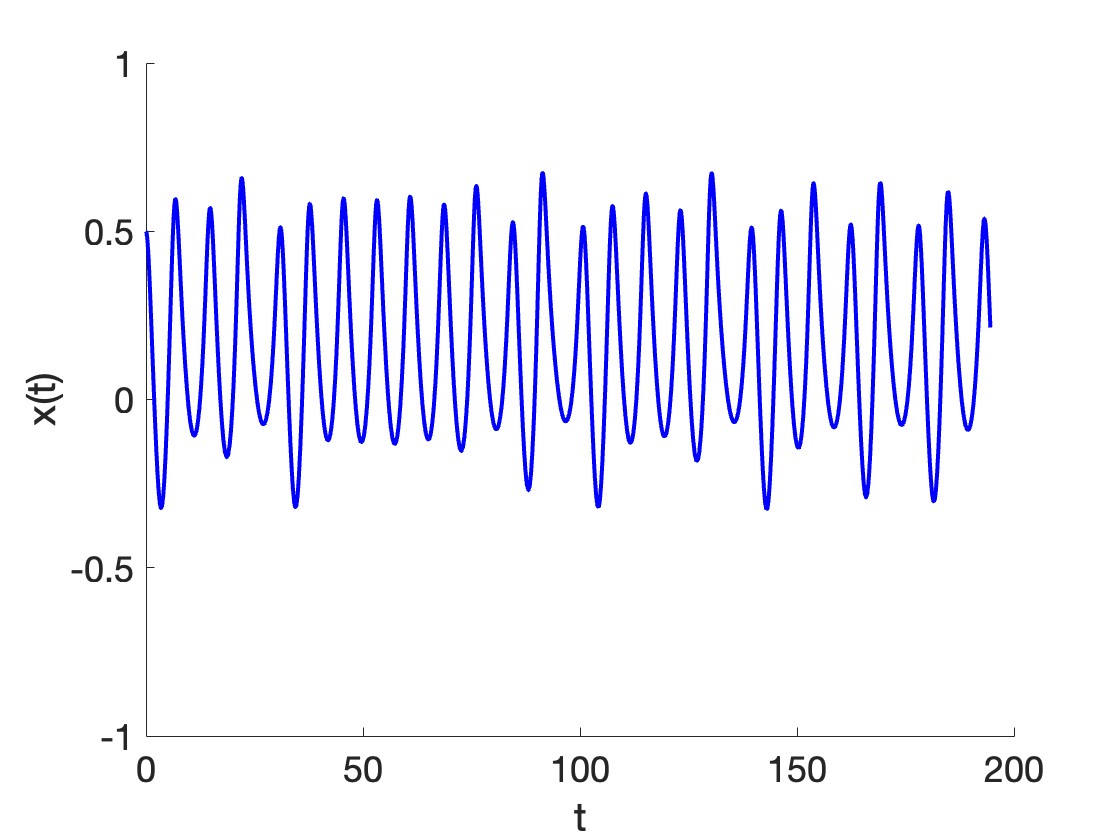}} 
\end{center}
\caption{National Product $Y(t)$ with parameter values 
$ s_r=0.1, \delta_k=0.15,  dr= 0.15,  \alpha = 0.2,  \beta=0.35$. 
The parameter $s_k=0.4$ produced in fig.~\ref{fig2c} shows economic growth without fluctuations of
$s_k$, it is reproduced here again tending to a constant value. In addition 2  oscillating chaotic 
results for $Y(t)$ are emerging for $s_k = 0.4 \pm \frac{1}{2} x(t)$. 
Right the chaotic timeseries $x(t)$. 
 \label{Ychaos}}
\end{figure} 
In an economy capital investments are not fixed but are changing daily as a consequence 
of economic and political changes. Assuming we want to model small fluctuations that take place 
when there are no exceptionally large changes in the economy, we modify system \eqref{sys1} 
by replacing capital investment coefficient $s_k$ by $s_k + c x(t)$.  As an example, the timeseries 
$x(t)$ is 
chosen of the chaotic system NE9 studied in \cite{BVchaos}. 
This system has the advantage that chaos is produced as the endproduct of a series of 
solutions with increasing 
periods, we find in the chaotic morion still some cyclic behaviour. In chaotic systems that 
arise at unstable global motion, like in the Lorenz attractor, the timeseries is less suitable for 
economic models. We modify system \eqref{sys1} to: 
\begin{eqnarray} \begin{cases} \label{sys1chaos} 
\frac{dK}{dt}&  = (s_k + cx(t)) E^{\alpha}K^{\beta} - \delta_k  K,\\
\frac{dE}{dt} & = s_r E^{\alpha}K^{\beta} - \delta_r E.
\end{cases} \end{eqnarray} 
For reference we present the equations governing the chaotic system NE9:
\begin{equation} \label{NE9}
\dot{x}=y, \dot{y} =-x -yz, \dot{z}= -xz +7x^2 -0.55.
\end{equation}  
We put initially $x(0)=0.5, y(0)=z(0)=0$. 
Replacing the number $0.55$ in system \eqref{NE9} by other positive numbers will produce 
very different dynamics.  
The average of $x(t)$ of the NE9 system over time zero 
to $t$ is of a generalised form as $x(t)$ is not periodic or quasi-periodic. We put for the average 
if $t >0$: 
\begin{equation} \label{average}
A(t) = \frac{1}{t} \int_0^t x(s)ds. 
\end{equation} 
We have that the average $A(100)$ is a
small positive number close to $0.14$, so $x(t)$ spends more time at positive values than at 
negative investment values. 
We compare 2 growth scenarios by choosing $c= \pm 0.5$. The plus sign describes a case where 
a hype cycle of new technical developments produces chaotically new investments. 
The minus sign refers 
to a case where erratic political measures entails uncertainty to invest. 

In fig.~\ref{Ychaos} we show growth to a stable National Product $Y$  by  the 
case $s_k=0.4 , c=0$.and the 2 cases with unpredictable, chaotic oscillations.  \\

\section{Controlling consumption} \label{sec3}

\begin{figure}[ht]  
\begin{center} \resizebox{!}{8cm}{\includegraphics{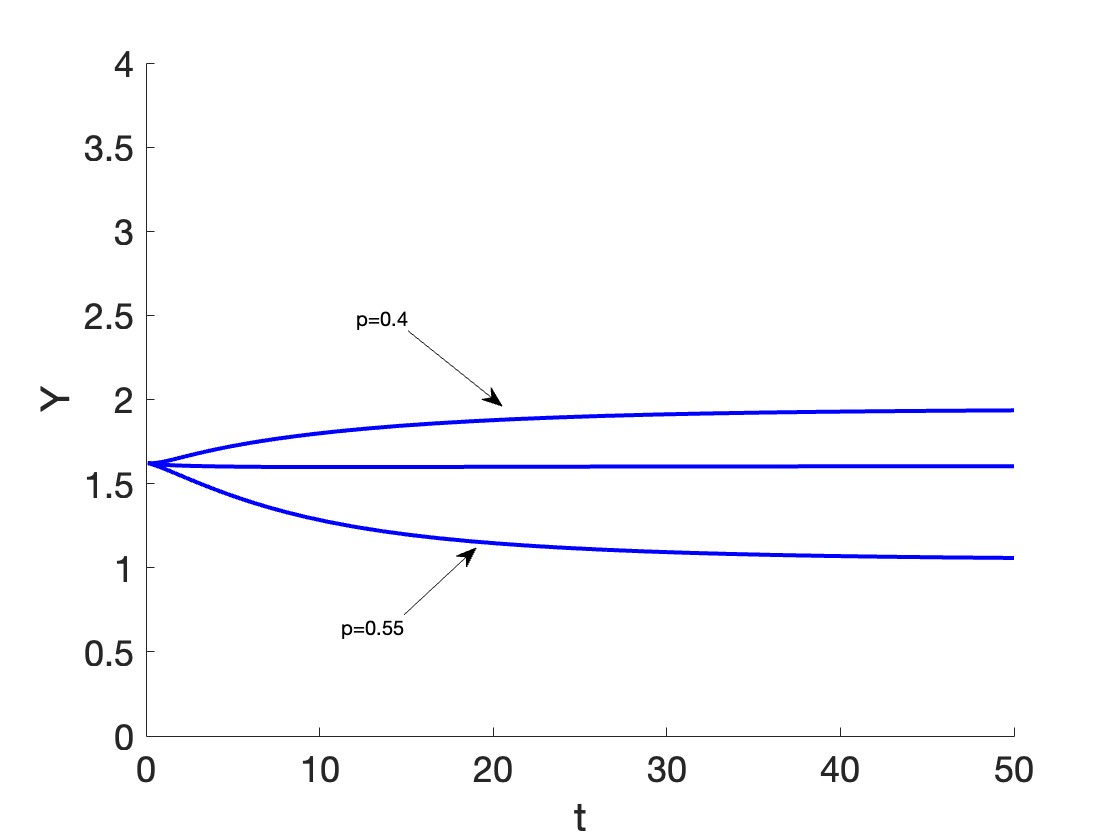}}
\end{center}
\caption{Dynamics of system \eqref{sys2} by controlling consumption. 
We have $E(0)=1, K(0)=4$ and parameters 
$s_k=0.4, \delta_k=0.15,  \delta_r=0.25, \alpha = 0.2,  \beta=0.35$. Investment in education and 
research starts at $s_r(0)=0.1$;
$s_r$ varies according to eq.~\eqref{sys2}. $Y(t)$ is shown for 3 rather different values of $p$, 
lowest $p=0.55$ with decrease of $Y(t)$, nearly constant growth for $p=0.47$ and highest growth 
at $p= 0.4$. With our choice of parameters $p= 0.47$ is a 
tipping point for consumption.
 \label{fig4a}}
\end{figure} 
We will investigate the consequences of changes of investments 
in education and research by government control. 

Already in 1961 E. Phelps \cite{P61} looked for conditions to obtain a fixed economic output 
to a fixed capital input ratio (the ``Golden Rule'').  An equilibrium solution can be obtained by 
suitable choices of 
all the parameters of a classical Cobb-Douglas function leading to so-called golden-age growth 
with ratio National Product and capital investments constant with time.  
In the ``Golden Rule''  model education and innovation investments are not seen as a separate type of 
investment, also stability questions have still to be solved at this stage. \\
As we will see, aiming at an optimal savings rate that maximizes per capita consumption in the 
long-run equilibrium 
is different from a direct control choosing a politically accepted consumption target as part of the 
increasing or decreasing National Product. \\

Suppose that in a new approach the government wants to guarantee consumption by controlling 
investment in education and research. There may be various reasons for this, for instance the 
financing of 
pensions  of an aging population, the increase of the expenses of medical care or defense or 
simply a negative populist view of education as an elite activity in general. 
So, for political reasons 
one aims then at consumption  $C$ as a certain fraction $p$ of the national 
product $Y$. A given target of consumption $C$ can be achieved by varying the investment 
factor $ s_r$ with time keeping capital investment factor $s_k$ constant. This results in the 
following equation to control investments: 
\begin{equation} \label{control}
 \frac{ds_r}{dt} = C- pY. 
 \end{equation} 
 
 \begin{figure}[ht]  
\begin{center}   \resizebox{!}{8cm}{
\includegraphics{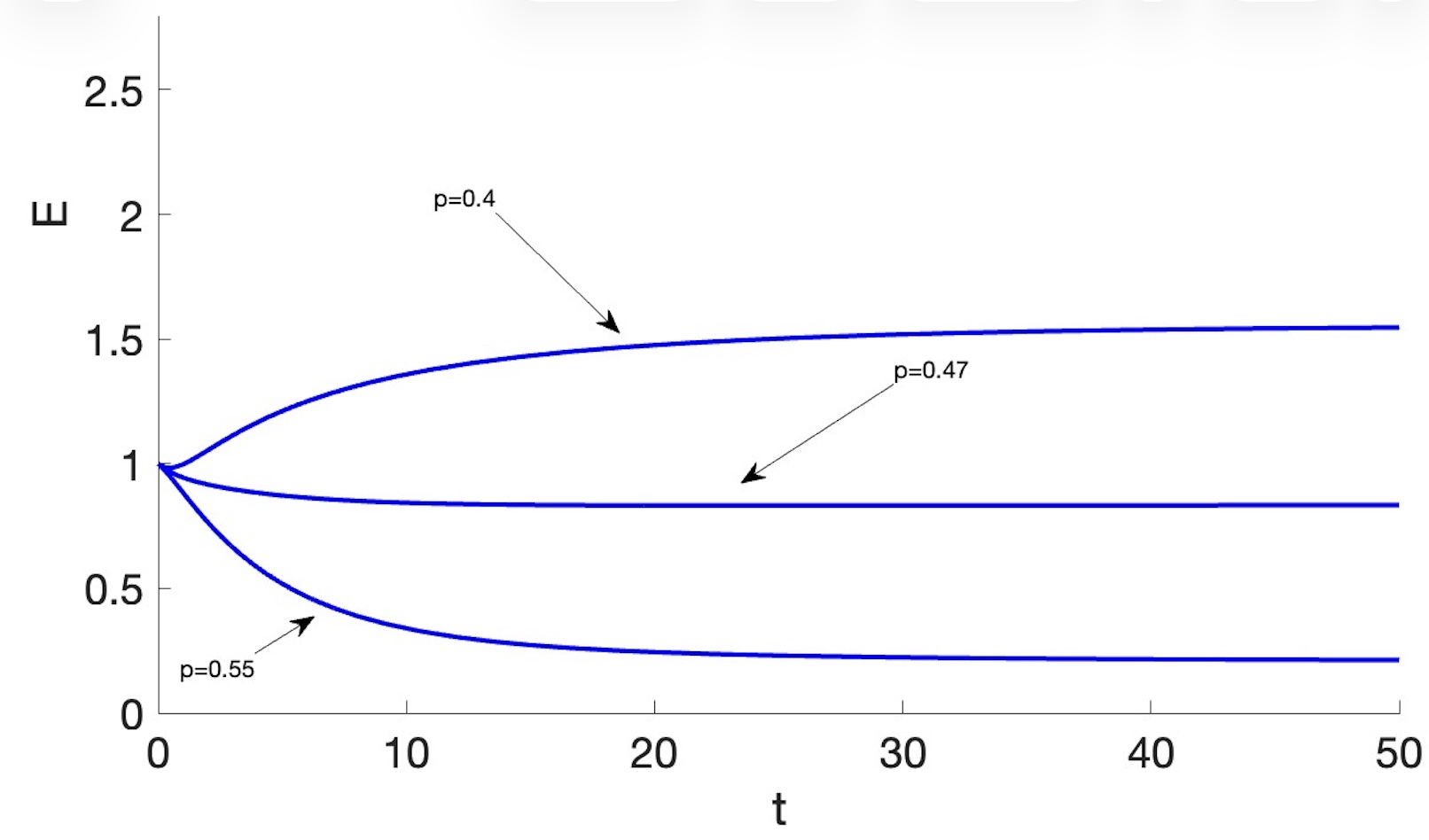}}
\end{center}
\caption{Dynamics of $E(t)$ in system \eqref{sys2} by controlling consumption. 
$E(t)$ is shown corresponding with the same 3 values of control $p= 0.4, 0.47. 0.55$ as 
in fig.~\ref{fig4a}. 
 \label{fig4b}}
\end{figure} 
If the consumption $C$ is larger than target $pY$, the investment coefficient $s_r$ for education 
and research 
has room to increase, if the consumption $C$ is smaller than target 
$pY$ then $s_r$ has to decrease.  This type of control is inspired by chemical physics, for 
references and other mathematical physics applications see \cite{VB24}.\\
 As $C= (1-s_k-s_r)Y$ we have with eqs \eqref{CB}-\eqref{control} and system \eqref{sys1} 
the  3-dimensional system: 
\begin{eqnarray} \label{sys2}
\begin{cases}
\frac{dK}{dt}&  = s_k E^{\alpha}K^{\beta} - \delta_k  K,\\
\frac{dE}{dt} & = s_r E^{\alpha}K^{\beta} - \delta_r E,\\
\frac{ds_r}{dt} & = (1-s_k-s_r-p)E^{\alpha}K^{\beta}. 
\end{cases}
\end{eqnarray} 

The equilibrium (critical point )given by \eqref{critpt} will change as $s_r$ is varying by  
changing parameter $p$ and the evolution of the economy; $E_0, K_0$ are again given by eq.~\eqref{critpt} and:
\begin{equation} \label{pcritsr} 
s_r = 1-s_k-p. 
\end{equation} 
Linearisation at the equilibrium produces the matrix:
\begin{equation} \label{pcrit}  \left( \begin{array}{ccc} 
(\beta  - 1) \delta_k & \alpha \frac{s_k}{s_r} \delta_r  & 0\\
\beta \frac{s_r}{s_k} \delta_k & (\alpha -1) \delta_r & E_0^{\alpha}K_0^{\beta} \\
0 & 0& - E_0^{\alpha}K_0^{\beta} 
\end{array}  \right) \end{equation} 
where $s_r$ satisfies eq.~\eqref{pcritsr}.
The characteristic equation has the same first 2 eigenvalues as for system \eqref{sys2} 
with $s_r$ given by eq.~\eqref{pcritsr}  but as 
a 3rd eigenvalue $\lambda_3 = -E_0^{\alpha}K_0^{\beta} $. So the controlled equilibrium is also 
stable. \\ 

In fig.~\ref{fig4a}  we show the growth of the National Product $Y(t)$ for the cases $p= 0.4, 0.47, 
0.55$. We show the evolution of $E(t)$ in the corresponding 3 cases of $p$ in fig.~\ref{fig4b}. 
We conclude that reducing the investment in education-research to increase 
the possibility of increase of consumption may affect economic growth negatively. A tipping point 
for the National Product with the parameters of fig.~\ref{fig4a} is close to $p= 0.47$. $Y(t)$ remains in this case near the starting value $1.62$ and consumption $C(t) \rightarrow 0.76$.
If $p=0.55$ we have that $Y(t) \rightarrow 1.2$ and $C(t) \rightarrow 0.66$ so trying a control 
aiming at higher consumption as part of the National Product may in fact reduce consumption. 
If $p=0.4$ we have that $Y(t) \rightarrow 1.95$ and $C(t) \rightarrow 0.78$.

\section{Discussion} 
\begin{enumerate} 
\item We have constructed a low-dimensional macro-economic model involving investment in 
education and research in section \ref{sec2}; this model is changed in section \ref{sec3} by 
adding a control to achieve consumption as a given fraction $p$ of National Income. The control 
is effective and has long time consequences for investments in education and research and 
subsequent increase or decrease of the National Product and consumption. 

\item It turns out in section \ref{sec2} that economic growth is also sensitive to the quality 
of education (learning something that matters and remaining well-informed) and innovation (high quality and up-to-date research). 
This result was obtained by small changes of the coefficient $\delta_r$ describing the persistence 
of useful knowledge and expertise while 
keeping all the other parameter values constant. 

\item
In the practice of capital investments there is not  one source of money but many fluctuating sources 
dependent on different circumstances. We can model this by introducing a chaotically fluctuating 
modulation of the investment factor $s_k$ with hype-induced investments and corresponding growth 
of National Product $Y$, also the case of erratic policies leading to decrease of investments. 
To study growth and decline of the economy it 
would be useful to have examples of realistic timeseries of such fluctuations.

\item An interesting result of section \ref{sec3} is that aiming at a lower fraction $pY$ of National 
Product $Y(t)$ and so relatively less consumption this can result in more growth and so in more 
real consumption.

\end{enumerate}

\subsection*{Acknowledgement}
The numerical results of this note were obtained by {\sc Matcont} ode 78 under {\sc Matlab}. 

\subsection*{Conflict of Interest}
The author has no conflicts of interest.

\subsection*{No funding}

\subsection*{Data Availability Statement}
The data that support the findings of this study are available from the corresponding author upon reasonable request.

\end{document}